\documentclass[a4paper,12pt]{amsart}

\usepackage{amsmath}
\usepackage{amsthm}
\usepackage{amsfonts}
\usepackage[cp1250]{inputenc}
\usepackage[T1]{fontenc}
\usepackage{epsfig}
\usepackage{setspace}

\newcommand{\CC}{{\mathbb C}}
\newcommand{\RR}{{\mathbb R}}

\newcommand{\ZZ}{{\mathbb Z}}

\newcommand{\ep}{\varepsilon}
\newcommand{\GG}{ G\mbox{-set}}
\newcommand{\GGs}{ G\mbox{-sets}}
\newcommand{\too}{\longrightarrow}
\newcommand{\TU}{T_{U}}
\newcommand{\TL}{T_{L}}
\newcommand{\TR}{T_{R}}
\newcommand{\TUL}{T_{UL}}
\newcommand{\TUR}{T_{UR}}
\newcommand{\TUG}{T_{U}(\Gamma)}
\newcommand{\TLG}{T_{L}(\Gamma)}
\newcommand{\TRG}{T_{R}(\Gamma)}
\newcommand{\TULG}{T_{UL}(\Gamma)}
\newcommand{\TURG}{T_{UR}(\Gamma)}
\newcommand{\ra}{\ensuremath{\frac{1}{r}(1,a,r-a)}}

\newcommand{\tTG}{\widetilde{\Theta}(\Gamma)}
\newcommand{\RRp}{\RR_{\ge 0}}
\newcommand{\Gyz}[1]{\Gamma^{\mathrm{yz}}_{#1}}
\newcommand{\Gx}{\Gamma^{\mathrm{x}}_{\phantom{0}}}
\newcommand{\wTheta}{\widetilde\Theta}

\DeclareMathOperator{\Hom}{Hom} \DeclareMathOperator{\SL}{SL}
\DeclareMathOperator{\GL}{GL} 
 \DeclareMathOperator{\Spec}{Spec}
 \DeclareMathOperator{\Hilb}{Hilb}
\DeclareMathOperator{\diag}{diag} \DeclareMathOperator{\V}{V}
\DeclareMathOperator{\wt}{wt}
\DeclareMathOperator{\wtG}{wt_{\Gamma}}
\DeclareMathOperator{\SSS}{S} 
\DeclareMathOperator{\GCD}{GCD} \DeclareMathOperator{\spann}{span}

\newcommand{\GH}{G\hbox{-}\Hilb \CC^3}\newcommand{\GHnic}{\Hilb^{G}}

\newtheorem{theorem}{Theorem}[section]
\newtheorem{lemma}[theorem]{Lemma}
\newtheorem{definition}[theorem]{Definition}
\theoremstyle{definition}
\newtheorem{remark}[theorem]{Remark}
\newtheorem*{xremark}{Remark}
\newtheorem{corollary}[theorem]{Corollary}

\newtheorem{notation}[theorem]{Notation}

\begin{document}

\title{The $G$-Hilbert scheme for $\frac{1}{r}(1,a,r-a)$}
\author{Oskar K\k{e}dzierski}
\address{Institute of Mathematics\\Warsaw University\\
 ul. Banacha 2\\
 02-097 Warszawa\\ Poland}

\email{oskar@mimuw.edu.pl} 

\thanks{Research supported by a grant of Polish MNiSzW (N N201 2653 33).}

\date{\today}
\maketitle

\begin{abstract}

Following Craw, Maclagan, Thomas and Nakamura's
work~\cite{Nakamura:Hilbert},\cite{CrawMaclaganII} on Hilbert
schemes for abelian groups, we give an explicit description of the
$\GHnic\CC^3$ scheme for
$G=\langle\diag(\ep,\ep^a,\ep^{r-a})\rangle$ by a classification of
all $G-$sets. We describe how the combinatorial properties of the
fan of $\GHnic\CC^3$  relates to the Euclidean algorithm.

\end{abstract}

%\doublespace

\section{Introduction}
For any finite, abelian subgroup $G$ of $\GL(n,\CC)$ of order $r,$
Nakamura defines the $G-$Hilbert scheme $\text{Hilb}^G \CC^n$ as the
irreducible component of the $G-$fixed set of the scheme
$\Hilb^{r}\CC^n$ which contains free orbits.

For such groups, the normalization of $\GHnic\CC^n$ is a toric variety.
The scheme $\GHnic\CC^n$ is described in~\cite{Nakamura:Hilbert} in terms of $G-$sets.
In fact, the description is carried by a classification of $G-$sets.

There are several known cases when $\GHnic\CC^n$ itself is a toric
variety (i.e.~it is normal): for $n=2$ and $G\subset\GL(2,\CC)$ by
Kidoh \cite{Kidoh:Hilbert}, for $n=3$ and $G\subset\SL(3,\CC)$ by
Craw and Reid \cite{CrawReid:GHSL3}, for any $n\ge 2$ and
$G=\langle \diag(\ep,\ep^2,\ep^4,\ldots,\ep^{2^n})\rangle$ by
Sebestean \cite{Sebestean:toric}. In all these cases, if  $n\ge 3$
the quotient $\CC^n /G$ has canonical, non--terminal singularities.

Craw, Maclagan and Thomas in~\cite{CrawMaclaganII}
describe $\GHnic\CC^n$  for any finite, abelian group
$G\subset\GL(n,\CC)$ in terms of initial ideals of some fixed
monomial ideal by varying weight order. This gives a numerical
method for finding the fan of $\GHnic\CC^n.$

In this paper, we use~\cite{Nakamura:Hilbert} and~\cite{CrawMaclaganII} to give a conceptual description of $\GHnic\CC^3$ scheme
for any cyclic subgroup $G\subset\GL(3,\CC)$ for which the quotient
$\CC^3 / G$ is a terminal singularity; this is Theorem~(\ref{algorithm}) below.
By Morrison and Stevens~\cite{MorrisonStevens},
any such group is conjugated to a group generated by a diagonal
matrix $\diag(\ep,\ep^a,\ep^{r-a}),$ where $a$ and $r$ are any
coprime natural numbers and $\ep$ is an $r-$th primitive root of
unity.

The description is carried out by classification of all possible
$G-$sets in families, called triangles of transformations. These
families correspond to steps in the Euclidean algorithm for $b$ and
$r-b,$ where $b$ is an inverse of $a$ modulo $r$ (see Main
Theorem~(\ref{algorithm})). We prove that there are
$\frac{1}{2}(3r+b(r-b)-1)$ different $\GGs$ (see Theorem~(\ref{the
number of G-sets})).

We show that for $a,r-a>1$ the $\GHnic\CC^3$ scheme is a normal
variety with quadratic singularities. Note that $\GHnic\CC^3$ for
$a=1$ or $r-a=1$ is isomorphic to the Danilov resolution of $\CC^3 /
G$ singularity by~\cite{Kedzierski:CohGH}.

The paper is organized as follows. Section~\ref{section:GHilb
basics} recalls basic definitions from~\cite{Nakamura:Hilbert}. Section~\ref{section:classification of
G-sets} contains classification of the $\GGs$ by the
number of valleys.  It is used to show that the $\GHnic$ is normal.
Section~\ref{section:G-isgaw transformations} contains definition
of a primitive $G-$set. Every such $G-$set gives rise to a family
of $G-$sets. The union of toric cones corresponding to $G-$sets in
such family is called a triangle of transformations.
In Section~\ref{triangles of transformation} we show how
to obtain a new primitive $G-$set from another one. In Sections~\ref{main
algorithm} the combinatoric
properties of primitive $G-$sets and the triangles of transformations are related to the Euclidean
algorithm.  We show that all subcones of cones in all triangles of transformations form the fan of $\GHnic$ scheme. The formula counting the number of $\GGs$ is
given at the end of Section~\ref{main algorithm}.
Section~\ref{section:example GHilbert} contains a concrete
example of $\GHnic$ scheme for $G\cong \ZZ_{14}.$

I would like to thank Professor Miles Reid for introducing me to
this subject.

\section{Basic definitions}\label{section:GHilb basics}

Let us fix two coprime integers $r,a\geq 2.$ Without loss of
generality we may assume that $a<r-a<r.$ Denote by $G$ the cyclic
group $\ZZ_r,$ considered as a subgroup of $\GL(3,\CC),$ generated
by matrix $\diag (\ep,\ep^a,\ep^{r-a}),$ where $\ep=e^{\frac {2\pi
i}{r}}.$  The group $G$ has $r$ characters which may be identified
with $1,\ep,\ep^2,\ldots,\ep^{r-1}.$

We follow the notation of~\cite{Nakamura:Hilbert}. Let $N_0=\ZZ e_1
\oplus \ZZ e_2 \oplus \ZZ e_3$ denote a free $\ZZ-$module with
$\ZZ-$basis $e_i.$ The lattice dual to $N_0$ will be denoted
$M_0=\Hom_\ZZ (N_0,\ZZ)=\ZZ e_1^* \oplus \ZZ e_2^* \oplus \ZZ e_3^*,
$ where $e_i^*(e_j)=\delta_{ij}.$ For the rest of this paper the
variables $x,y,z$ will be identified with $e_1^*,e_2^*,e_3^*$ and a
multiplicative notation will be used in the lattice $M_0.$ For
example, vector $2e_1^*-e_3^*$ will be identified with the Laurent
monomial $x^2z^{-1}.$

Let $M_0^0$ be the positive octant in $M_0,$ identified with
monomials in the ring $\CC[x,y,z].$ Set $N=N_0+\ZZ
\frac{1}{r}(e_1+ae_2+(r-a)e_3)$ and let $M=\Hom_\ZZ (N,\ZZ)$ be a
dual lattice. Lattice $M$ will be identified with a sublattice of
$M_0$ consisting of $G-$invariant Laurent monomials. When no confusion
arise, vector $a_1 e_1 + a_2 e_2 + a_3 e_3$ will be denoted
$(a_1,a_2,a_3).$ For example $\frac{1}{5}(1,2,3)$ stands for
$\frac{1}{5}e_1+\frac{2}{5}e_2+\frac{3}{5}e_3.$

Let $G^{\vee}$ denote the character group of $G.$ The group $G$ acts on
the left on regular functions on $\CC^3$ by setting $(g\cdot
f)(p)=f(g^{-1}p),$ where $g\in G,$ $p\in\CC^3$ and $f$ is a regular
function on $\CC^3.$ This action can be extended to the lattice
$M_0$ (by identifying $M_0$ with the lattice of exponents of Laurent monomials in $x,y,z$).
Thus, we have the natural grading:
\[M_0=\bigoplus_{\chi \in G^{\vee}}M^{\chi}_0.\]

\begin{definition}
Let $\wt:M_0 \too G^{\vee}$ denote group homomorphism sending an element of
the lattice $M_0$ to its grade.
\end{definition}

We will denote by $m \mod n$ an integer $k\in{0,\ldots,n-1}$ such
that $n|(m-k).$

\begin{definition}(Nakamura)\label{defintion of G-set}
A subset $\Gamma$ of monomials in $\CC [x,y,z]$ is called a
\emph{$G-$set} if
\begin{enumerate}
\item{it contains the constant monomial $1,$}
\item{if $vw \in \Gamma $ then $v \in \Gamma $ and $w \in \Gamma ,$}
\item{the restriction of the function $\wt$ to $\Gamma$ is a bijection.}
\end{enumerate}
\end{definition}

\begin{xremark} \label{remarks:possible clusters}
Since $\wt(1)=\wt(yz),$ it follows that $yz\notin\Gamma$ for any
$\GG$ $\Gamma.$ Hence the monomials in $\Gamma$ are of the form
$x^*y^*$ and $x^*z^*,$ where $*$ stands for any nonnegative integer.
\end{xremark}

\begin{definition}\label{definition:ijk}
For any $\GG$ $\Gamma$ define $i(\Gamma),j(\Gamma),k(\Gamma)$ to be
the unique nonnegative integers such that
$$x^{i(\Gamma)}\in\Gamma,x^{i(\Gamma)+1}\notin\Gamma,$$
$$y^{j(\Gamma)}\in\Gamma, y^{j(\Gamma)+1}\notin\Gamma,$$
$$z^{k(\Gamma)}\in\Gamma, z^{k(\Gamma)+1}\notin\Gamma.$$
When no confusion arise we write for short:
$$i=i(\Gamma),$$
$$j=j(\Gamma),$$
$$k=k(\Gamma).$$
\end{definition}

%On pictures a $\GG$ is presented as a plane lattice of boxes with
%monomial in every box such that the exponents at $x, y,z$ increase
%in the upper, right, left direction, respectively (cf.
%Figure~\ref{figure:sample G-set}).

\begin{definition}(Nakamura)\label{definition:definition of valley}
A monomial $x^m y^n$ (resp. $x^m z^n$) for $m,n\ge 0$ is called a
\emph{$y-$valley} (resp. \emph{$z-$valley}) for $\Gamma,$ if
$$x^m y^n,x^{m+1} y^n,x^m y^{n+1}\in\Gamma
\quad\mbox{but}\quad x^{m+1} y^{n+1}\notin\Gamma$$
$$\text{(resp. } x^m z^n,x^{m+1} z^n,x^m z^{n+1}\in\Gamma\ \text{ but }\ x^{m+1}
z^{n+1}\notin\Gamma).$$ We call a $y-$valley or $z-$valley a valley
for brevity.
\end{definition}

%\begin{figure}[h] \centering \includegraphics[width=0.1\textwidth]{123Gset.eps}
%\caption{An example of a $\GG$ for $r=5,a=2.$} \label{figure:sample
%G-set}
%\end{figure}

\begin{definition}

For any $v\in M_0^0$ let $\wtG (v)$ denote the unique $w\in
\Gamma$ such that $\wt(v)=\wt(w).$
\end{definition}

\section{Classification of $\GGs$}\label{section:classification of
G-sets}

In this section we show that any $\GG$ has at most one $y-$valley
and at most one $z-$valley. Following Nakamura, for every $G-$set we
construct a semigroup $\SSS(\Gamma)$ in the lattice $M$ and prove
that it is saturated.  It turns out that the $\GGs$ correspond to
the cones of maximal dimension in the fan of $\Hilb^G {\CC^3}.$

\begin{remark}\label{remark:tricks}
The following statements are immediate from the definitions:
\begin{enumerate}
\item{if $\wtG(v)=w,\ v\notin\Gamma$ and $u \cdot w\in \Gamma,$ then
$u \cdot v\notin \Gamma,$}
\item{if $\wtG(v)=w,$ then $\wtG(u\cdot v)=u\cdot w$ for any $u\in M_0$ such that $u\cdot w\in \Gamma,$}
\item{if $\wtG(v)=w,\ u\in M$ then $\wtG(u\cdot v)=w.$}
\end{enumerate}
\end{remark}

\begin{corollary}\label{C:extremes}
Let $\Gamma$ be a $\GG$ and $v\in M_0^0-\Gamma.$ If $x^{-1}\cdot
v\in\Gamma$ (resp. $y^{-1}\cdot v\in\Gamma,z^{-1}\cdot v\in\Gamma$)
then $\wtG(v)=w,$ where $w\in\Gamma$ but $x^{-1}\cdot w\notin\Gamma$
(resp. $z\cdot w\notin\Gamma,y\cdot w\notin\Gamma$).
\end{corollary}

\begin{proof}
Use observation (1) and (3) from Remark~(\ref{remark:tricks}).
\end{proof}

\begin{lemma}\label{number of valleys}

A $\GG$ can only have $0,1$ or $2$ valleys.

\end{lemma}

\begin{proof}
Suppose that $x^m y^n$ is a $y-$valley for $\Gamma.$ Then
$v=x^{m+1}y^{n+1}$ satisfies assumptions of
Corollary~(\ref{C:extremes}). Hence, $x^{-1}\cdot \wtG(v)
\notin\Gamma$ and $z\cdot \wtG(v)\notin \Gamma,$ so
$\wtG(v)=z^{k(\Gamma)}.$ Therefore, $\GG$ $\Gamma$ has at most one
$y-$valley, and, analogously at most one $z-$valley.
\end{proof}

\begin{corollary}\label{where extremes go}
Suppose that $\GG$ $\Gamma$ has $y-$valley $w$ and $z-$valley $v.$
Then
$$\wtG(y^{j(\Gamma)+1})=x\cdot w,$$
$$\wtG(z^{k(\Gamma)+1})=x\cdot v.$$

\end{corollary}

\begin{proof}
Use observation (2) from Remark~(\ref{remark:tricks}).
\end{proof}

\begin{notation}
From now on we will usually denote by $i_y,j_y$ the exponents of the
$y-$valley $x^{i_y}y^{j_y}$ and by $i_z,k_z$ the exponents of the
$z-$valley $x^{i_z}z^{k_z}$ of some fixed $\GG$ $\Gamma.$
\end{notation}

\begin{lemma}\label{G-sets with no valleys}
The only possible $\GGs$ with no valleys are:
$$\Gx =\{1,x,\ldots,x^{r-1}\},$$
$$\Gyz{l} =\{y^{r-l-1},\ldots,y,1,z,\ldots,z^{l}\}\ \mbox{for}\ l=0,\ldots r-1 .$$
\end{lemma}

\begin{proof}

Let $i,j,k$ be integers like in Definition~(\ref{definition:ijk}).
Corollary~(\ref{C:extremes}) shows that $\wtG (y^{j+1})=x^{i'}z^k,$
for some $i'\ge 0.$ If $i'=0,$ then $\wt(z^{k+1})=\wt(y^j)$ and
since $a,r$ are coprime, it follows that $j=r-k-1,$ hence $i=0.$
Consider the case $i'>0.$ Then $\wtG(x^{i'-1}z^{k+1})=x^{i''} y^j$
by Corollary~(\ref{C:extremes}). It follows immediately that
$i''=i=r-1$ and so $j=k=0.$
\end{proof}

\begin{lemma}\label{extremes in case of 1 valley}
Let $\Gamma$ be a $\GG$ with exactly one valley. If $\Gamma$ has
$y-$valley equal to $x^{i_z}z^{k_z},$ then
$$\wtG(x^{i+1})=z^{k-k_z},$$
$$\wtG(z^{k+1})=x^{i-i_z}y^j.$$
If $\GG$ $\Gamma$ $z-$valley equal to $x^{i_y}y^{j_y},$ then
$$\wtG(x^{i+1})=y^{j-j_y},$$
$$\wtG(y^{j+1})=x^{i-i_y}z^k.$$
\end{lemma}

\begin{proof}
We prove the lemma in the case of $z-$valley $w=x^{i_z}z^{k_z}.$ The
monomial $\wtG(z^{k+1})$ is of the form $x^l y^j,$ where $0\le l \le
i.$ Noting that $\wtG(xz\cdot w)=y^j$ we get $l=i-i_z.$ It follows
that the monomials $x^{i-i_z}y^j$ and $z^{k+1}$ are of the same
weight, therefore $\wtG(z^{k+1})=x^{i-i_z}y^j.$
\end{proof}

\begin{lemma}\label{extremes in case of 2 valleys}
Let $\Gamma$ be a $\GG$ with two valleys $v,w,$ where
$$v=x^{i_y}y^{j_y},$$
$$w=x^{i_z}z^{k_z}.$$
Then $i_y+i_z+1=i,$ and
\begin{displaymath}
 \wtG(x^{i+1}) = \left\{
 \begin{array}{ll}
y^{(j-j_y)-(k-k_z)} &\; \text{if}\ \ (j-j_y)-(k-k_z)\geq 0,\\
z^{(k-k_z)-(j-j_y)} &\; \text{otherwise.} \\
 \end{array} \right.
\end{displaymath}
\end{lemma}

\begin{proof}
Let $u$ be a monomial such that $u\notin\Gamma$ and
$x^{-1}u\in\Gamma.$ Then $\wtG(u)=z^l$ for some $0 \le l\le k$ or
$\wtG(u)=y^l$ for $0 \le l\le j.$ We know already that $\wtG(xz\cdot
w)=y^j$ and $\wtG(xy\cdot v)=z^k,$ which implies that $\wtG(x^{i+1})
= y^{(j-j_y)-(k-k_z)}$ if $(j-j_y)-(k-k_z)\geq 0$ and
$\wtG(x^{i+1})=z^{(k-k_z)-(j-j_y)}$ otherwise. The monomial
$x^{i_y+i_z+1}$ has the same weight as $x^{i+1}$ hence they are
equal.
\end{proof}

\begin{definition}(Nakamura) \label{cones for G-set}
For any $v\in M_0$ and a $\GG$ $\Gamma$ define (using a
multiplicative notation in the lattice $M_0$)
\[s_\Gamma(v)=v\wt_\Gamma^{-1}(v).\]
We will write it simply $s(v)$ when no confusion can arise.
 Define the cones
\[\sigma (\Gamma) = \{ \alpha\in N_0 \otimes_{\ZZ} \RR\; |\; \langle
\alpha,s_\Gamma(v)\rangle\ge 0,\quad \forall v\in M_0^0 \},\]
\[\sigma^{\vee}(\Gamma)=\{ v\in M_0 \otimes_{\ZZ} \RR\; |\;
\langle\alpha,v\rangle \ge 0,\quad \forall
\alpha\in\sigma(\Gamma)\},\] where $\langle\cdot,\cdot\rangle$
denotes the pairing between $N_0$ and $M_0.$

Let $\SSS(\Gamma)$ be a subsemigroup of the lattice $M,$ generated
by the set $\{s_\Gamma(v)\in M\; |\; v\in M_0^0\}$ as a semigroup.
Set
$$\V(\Gamma)=\Spec \CC[\SSS(\Gamma)].$$
\end{definition}
Note that
$$\CC[\SSS(\Gamma)]\subset\CC[\sigma^{\vee}(\Gamma)\cap M].$$
 Moreover, the cones
$\sigma(\Gamma),\sigma^{\vee}(\Gamma)$ are dual to each other and
the cone $\sigma^{\vee}(\Gamma)\cap M$ is the saturation of the
semigroup $\SSS(\Gamma)$ in the lattice $M.$ It will follow from
Lemma~(\ref{Nakamura's lemma}) that $S(\Gamma)$ is finitely
generated as a semigroup.

%\begin{example}
%Let $\Gamma$ be the $\GG$ from Figure~\ref{figure:sample G-set}. The
%cone $\sigma^{\vee}(\Gamma)$ is generated by monomials
%${x^2}y^{-1},x^{-1}y^2z^{-1},x^{-1}{z^2}$ and the cone
%$\sigma(\Gamma)$ is generated by $\frac{1}{5}(1,2,3),$
%$\frac{1}{5}(2,4,1),$ $\frac{1}{5}(4,3,2)$ over $\RR_{\ge 0}.$  The
%primitive vectors along rays of $\sigma(\Gamma)$ form a $\ZZ-$basis
%of the lattice $N,$ thus $\V(\Gamma)\cong\CC^3.$
%\end{example}

\begin{theorem}[Nakamura]\label{ZNakamury}
Let $G$ be a finite abelian subgroup of $\GL(3,\CC).$  When $\Gamma$
varies through all $\GGs$ the set of all faces of all $3-$dimensional cones $\sigma(\Gamma)$
forms a fan in lattice $N\otimes\RR$ supported on the positive octant. Toric variety defined by this fan is
isomorphic to the normalization of the $\Hilb^G \CC^3$ scheme
%a normalization of the coherent component of $\GH$ scheme
(see~\cite[Theorem~2.11]{Nakamura:Hilbert} and \cite[$\S
5$]{CrawMaclaganI}). Moreover, the affine varieties $\V(\Gamma)$ form an
open covering of the $\Hilb^G \CC^3$ scheme when $\Gamma$ varies
through all $\GGs.$
\end{theorem}

\begin{lemma}(Nakamura)\label{Nakamura's lemma}
Let $A\subset M_0^0- \Gamma$ be a finite set such that
$M_0^0-\Gamma=A\cdot M_0^0.$ If $\sigma(\Gamma)$ is a
$3-$dimensional cone then $\SSS(\Gamma)$ is generated by the finite
set $\{\;s_\Gamma(v)\mid v\in A\;\}$ as a semigroup
(see~\cite[Lemma~1.8]{Nakamura:Hilbert}).
\end{lemma}

\begin{remark}
Note that Theorem~(\ref{ZNakamury}) and Lemma~(\ref{Nakamura's
lemma}) are stated in~\cite{Nakamura:Hilbert} without the assumption
on dimension of $\sigma(\Gamma)$ in which case they are false. A
counterexample and a correction can be found in~\cite[Example~4.12
and Theorem~5.2]{CrawMaclaganII}.
\end{remark}

\begin{lemma}\label{G-Hilb is normal}
Suppose that $\Gamma$ is a $\GG$ in the case of $\ra$ action. Then
the cone $\sigma(\Gamma)$ is $3-$dimensional. Moreover, if $\Gamma$
has $0$ or $1$ valley then %$\SSS(\Gamma)$ has $3$ generators and
$\SSS(\Gamma)\cong\CC[x,y,z].$ If $\Gamma$ has $2$ valleys then
%$\SSS(\Gamma)$ has $4$ generators and
$\SSS(\Gamma)\cong  \CC
[x,y,z,w]/(xy-zw).$ \
\end{lemma}

\begin{proof}
The lemma will be proven only in the case of a $\GG$ with $2$
valleys as the method carries over to the other cases.

Suppose that $\Gamma$ is a $\GG$ with $2$ valleys, $v=x^{i_y}
y^{j_y},\ w=x^{i_z}z^{k_z}$ and set
\begin{align*}
\alpha&=x^{i+1},\\
\beta&=y^{j+1},\\
\gamma&=z^{k+1},\\
\delta_y&=xy\cdot v,\\
\delta_z&=xz\cdot w,
\end{align*}
where $i,j,k$ are the largest exponents such that $x^i,y^j,z^k$
belong to $\Gamma.$ We will start by showing that
$s(\beta),s(\gamma),s(\delta_y)$ and $s(\delta_z)$ generate
semigroup $\SSS(\Gamma).$ Assume that $u\in M_0^0,\; t=x,y$ or $z$
and note that
\[s(t\cdot u)=s(u)s(t\cdot\wtG(u)).\]
By the above formula it suffices to show that for any $u\in\Gamma$
such that $t\cdot u\notin \Gamma$ the Laurent monomial $s(t\cdot u)$
can be expressed as a product of $s(\beta),s(\gamma),s(\delta_y)$
and $s(\delta_z)$ with nonnegative exponents. By
Lemma~(\ref{extremes in case of 2 valleys}):
\[s(\alpha)= \left\{
 \begin{array}{ll}
x^{i+1}{y^{-(j-j_y)+(k-k_z)}} &\; \text{if}\ \ (j-j_y)\geq(k-k_z),\\
x^{i+1}{z^{(j-j_y)-(k-k_z)}} &\; \text{otherwise,} \\
 \end{array} \right. \]
$$s(\beta)={xy^{-(j+1)}\cdot w},$$
$$s(\gamma)={xz^{-(k+1)}\cdot v},$$
$$s(\delta_y)={xyz^{-k}\cdot v},$$
$$s(\delta_z)={xy^{-j}z\cdot w},$$
hence
\[s(\beta)s(\delta_z)=s(\gamma)s(\delta_y)=s(yz),\]
\[s(\alpha)= \left\{
 \begin{array}{ll}
s(\delta_y)s(\delta_z)(yz)^{j-j_y-1} &\; \text{if}\ \ (j-j_y)\geq(k-k_z),\\
s(\delta_y)s(\delta_z)(yz)^{k-k_z-1} &\; \text{otherwise.} \\
 \end{array} \right. \]
%(note that $j-j_y-1\geq 0$ and $k-k_z-1\geq 0$).

Let $u\in\Gamma$ and $y\cdot u \notin\Gamma.$ If $u=x^l y^j,$ where
$l=0,\ldots,i_y$ then $s(y\cdot u)=s(\beta).$ If $u=x^l y^j_y,$
where $l=i_y+1,\ldots,i$ then $s(y\cdot u)=s(\delta_y).$ Analogously
$s(z\cdot u)$ is equal to $s(\gamma)$ or  to $s(\delta_y)$ for any
$u\in\Gamma,z\cdot u\notin\Gamma.$

It remains to consider $u\in\Gamma$ such that $x\cdot
u\notin\Gamma.$ Observe that $\wtG(x\cdot u)$ is of the form $y^l$
or $z^l$ for some positive $l$ ($l=0$ can happen only if
$\Gamma=\Gx$). If $u'=y^{-1}u\in \Gamma$ then $x\cdot
u'\notin\Gamma$ and \[s(x\cdot u)=s(y\cdot xu')=s(xu')s(y\wtG(x\cdot
u'))=s(xu')(yz)^n,\text{ where }n=0,1.\] By induction for any such
$u\in\Gamma$ the monomial $s(x\cdot u)$ is equal to $p\cdot(xy)^m,$
where $m>0$ and $p=s(\alpha),s(\delta_y)$ or $s(\delta_z).$

This shows that $\SSS(\Gamma)$ is generated by
$s(\beta),s(\gamma),s(\delta_y)$ and $s(\delta_z).$ To conclude it
is enough to show that some (in fact any) $3$ out of $4$ generators
form a $\ZZ-$basis of the lattice $M.$ This is implied by computing
the following determinant, using equality from Lemma~(\ref{extremes
in case of 2 valleys}):

$$ \left|
                 \begin{array}{ccc}
                   -i_z-1 & j+1 & -k_z \\
                   i_y+1 & j_y+1 & -k \\
                   -i_y-1 & -j_y & k+1 \\
                 \end{array}
\right|=r.$$

\end{proof}

\begin{corollary}\label{G-Hilb is normal 2}
The semigroup $\SSS(\Gamma)$ coincides with the semigroup algebra
$\CC[\sigma^{\vee}(\Gamma')\cap M]$ for any $\GG.$ In particular, $\GHnic \CC^3$  is normal.
\end{corollary}

\section{$G-$igsaw transformations}\label{section:G-isgaw
transformations}

To get an effective description of the fan of the $\GHnic$
scheme, we introduce Nakamura's $G-$igsaw transformation,
 which will allow to organize $\GGs$ in families and to explain how
these are related to each other.

$G-$igsaw transformation is a method of constructing a new $\GG$ from the other.
In fact, two $\GGs$ $\Gamma$ and $\Gamma'$ are related by a $G-$igsaw
transformation if and only if the cones $\sigma(\Gamma)$ and $\sigma(\Gamma')$
share a $2-$dimensional face.

When reading Sections~\ref{section:G-isgaw transformations}
through~\ref{main algorithm}, it may be useful for a reader to consult an example provided in
Section~\ref{section:example GHilbert}.

\begin{lemma}(Nakamura)\label{lemma:existence of monom}
Let $\Gamma$ be a $\GG$ for the action of type $\frac{1}{r}(1,a,r-a)$ and let $\tau$ be a $2-$dimensional face of $\sigma(\Gamma).$ There exist two monomials $u\in M_0^0$ and $v\in\Gamma$ such that
\begin{enumerate}
\item $v=\wtG(u),$
\item $u,v$ do not have common factors in $M_0^0,$
\item $uv^{-1}$ is a primitive monomial,
\item $\tau=\sigma(\Gamma)\cap (uv^{-1})^{\bot},$
\end{enumerate}

\end{lemma}

\begin{proof}
This is a particular case of~\cite[Lemma~2.5]{Nakamura:Hilbert}
\end{proof}

\begin{definition}(Nakamura)\label{definition of Gigsaw}
Let $\Gamma$ be a $\GG$ and let $\tau$ be a $2-$dimensional face of
$\sigma(\Gamma).$ Suppose that monomials $u,v$ given by
Lemma~(\ref{lemma:existence of monom}) are not equal to $1$ and set
$c(w)=\max \{c\in\ZZ \mid {w}{v^{-c}} \in M_0^0 \}$ for any
$w\in\Gamma.$ We define the \emph{$G-$igsaw transformation of
$\;\Gamma$} in the direction of $\tau$ to be the set
$$ \Gamma'=\{w\cdot u^{c(w)}v^{-c(w)}\mid w\in\Gamma  \}.$$
\end{definition}

\begin{lemma}(Nakamura)
The $G-$igsaw transformation of a $\GG$ is a $\GG$.
\end{lemma}

\begin{proof}
See~\cite[Lemma~2.8]{Nakamura:Hilbert}
\end{proof}

\begin{lemma}\label{directions of G-isgaw transformations}
Suppose  that $\Gamma$ is a $\GG$ for the action
$\frac{1}{r}(1,a,r-a).$ Let
$\alpha=x^{i+1},\beta=y^{j+1},\gamma=z^{k+1},$ where $i,j,k$ are the
maximal exponents such that $x^i,y^j,z^k\in\Gamma.$ Let $\tau$ be a
$2-$dimensional face of $\sigma(\Gamma)$ and let $u$ be the monomial
given by Lemma~(\ref{lemma:existence of monom}). If $\Gamma$ has $0$
or $1$ valley then $u=\alpha,\beta$ or $\gamma.$ If $\Gamma$ has $2$
valleys then $u=\beta,\gamma,\delta_y$ or $\delta_z,$ where
$\delta_y$ is equal to the $y-$valley of $\Gamma$ multiplied by $xy$
and $\delta_z$ is equal to the $z-$valley of $\Gamma$ multiplied by
$xz.$

\end{lemma}

\begin{proof}
Suppose that $\Gamma$ has one valley and $\tau$ is a face of
$\sigma(\Gamma)$ dual to the ray of $\sigma^{\vee}(\Gamma)$ spanned
by $s(\alpha).$ The $1-$dimensional lattice $M\cap \tau^\bot$ has
$2$ generators. Therefore $uv^{-1}$ is equal either to $s(\alpha)$
or $s(\alpha)^{-1}.$ Clearly, the only choice is $u=\alpha,$
$v=\wtG(\alpha).$ Suppose that $d\in M_0^0$ is a common factor of
$u$ and $v.$ Then both $ud^{-1},vd^{-1}$ belong to $\Gamma$ and they
are of the same weight. Hence $d=1.$
\end{proof}

\begin{definition}
Let $\Gamma$ be a $\GG$ with $0$ or $1$ valley and let $\tau$ be the
$2-$dimensional face of $\sigma(\Gamma).$ The $G-$igsaw
transformation of $\Gamma$ in the direction of $\tau$ is called
\emph{upper (resp. right, left) transformation} if $u=\alpha$ (resp.
$u=\beta,u=\gamma.$), where the monomial $u$ is as in
Lemma~(\ref{lemma:existence of monom}). The upper, left and right
transformations of $\Gamma$ will be denoted by $\TUG,\TRG$ and
$\TLG,$ respectively.

By slight abuse of notation, the $G-$igsaw transformation of $\GG$
$\Gamma$ with $2$ valleys is called \emph{left (resp. upper left,
right, left) transformation} if the corresponding monomial $u$ is
equal to $\beta$ (resp. $\gamma,\delta_y,\delta_z$). The right,
left, upper right and upper left $G-$igsaw transformations of
$\Gamma$ will be denoted by $\TURG,\TULG,\TRG,\TLG,$ respectively.
%(see Figure~\ref{figure:jigsaw}).

\begin{definition}
We say that a $\GG$ $\Gamma$ is \emph{spanned} by monomials
$u_1,\ldots,u_n$ if $\Gamma$ consists of all monomials dividing
$u_1,\ldots,u_n.$ If $\GG$ $\Gamma$ is spanned by monomials
$u_1,\ldots,u_n$ we write
\[\Gamma=\spann (u_1,\ldots,u_n).\]
\end{definition}

\end{definition}

%\begin{example}
%For $r=14$ and $a=5$ the $\GG$
%$$\Gamma=\spann (x^3z^2,z^4)$$
%has one $z-$valley equal to $z^2$ (see Figure~\ref{figure:one valey
%transformation}). The upper transformation of $\Gamma$ are obtained
%by setting $u=x^4$ and $v=z^2$ in the definition of $G-$igsaw
%transformation. Thus, the upper transformation of $\Gamma$ is
%obtained by replacing each monomial $w\in\Gamma,$ divisible by
%$(z^{2n}$ but not by $(z^{2(n+1)},$  by the monomial
%$x^{4n}z^{-2n}\cdot w.$
%\end{example}

%\begin{figure}[h] \centering \includegraphics[width=0.82\textwidth]{new1valley_poprawiony.eps}
%\caption{$G$-igsaw transformations of a $\GG$ for
%$\frac{1}{14}(1,5,9).$} \label{figure:one valey transformation}
%\end{figure}

\begin{lemma}\label{classification of upper transformations}
Let $\Gamma=\spann (x^{i_y}y^j,x^i z^k),$ where $i_y<i$ (resp. let
$\Gamma=\spann (x^i y^j,x^{i_k}z^k),$ where $i_z<i$) be a $\GG$ with
one $y-$valley equal to $x^{i_y}$ (resp. one $z-$valley equal to
$x^{i_z}$).

Then \begin{align*}\TUG&=\spann(x^{i+i_y+1},x^{i_y}y^{j-1},x^i
z^k)\\
 \text{(resp. }\TUG&=\spann(x^{i+i_z+1},x^{i}y^{j},x^{i_k} z^{k-1})\text{)}.
\end{align*}
In particular,the upper transformation of $\Gamma$ has

\begin{itemize}
\item no valleys if and only if $\ j=1,k=0\ $ (resp. $j=0,k=1$). In
fact, in this case $\TU(\Gamma)=\Gx.$

\item one $z-$valley (resp. one $y-$valley) if and only if $\
j=1,k>0\ $ (resp. $j>0,k=1$). In both cases the valley is equal to
$x^i.$

\item two valleys: the $y-$valley equal to $x^{i_y}$ and the $z-$valley equal to
$x^{i}$ (resp. the $y-$valley equal to $x^{i}$ and the $z-$valley
equal to $x^{i_z}$) in the remaining cases.
\end{itemize}
\end{lemma}

\begin{proof}
The upper transformation is obtained by replacing each monomial $w\in\Gamma,$
divisible by $y^j$ (resp. by $z^k$) by the monomial $x^{n(i+1)}y^{-nj}\cdot w$
for some $n\ge 1.$ The proof is straightforward.
\end{proof}

\begin{lemma}\label{transformations of 2 valley into 2 valley}
Let $\Gamma$ be a $\GG$ with 2 valleys: $y-$valley equal to
$v=x^{i_y}y^{j_y}$ and $z-$valley equal to $w=x^{i_z}z^{k_z}.$
Assume that $\Gamma$ is spanned by $x^i y^{j_y},x^i z^{k_z},x^{i_y}y^j,x^{i_z}z^k.$
Let $T$ stand for right, left, upper right or upper left
transformation.

Then $T(\Gamma)$ is spanned by:
\[\begin{array}{llllcl}
  x^{i} y^{j_y}, & x^{i} z^{k_z-1}, & x^{i_y}y^{j+1}, & x^{i_z}z^{k} & \quad & T=\TR,k_z\ge 1\\
  x^{i} y^{j_y-1}, & x^{i} z^{k_z}, & x^{i_y}y^{j}, & x^{i_z}z^{k+1} & \quad & T=\TL,j_y\ge 1\\
   x^{i} y^{j_y+1}, & x^{i} z^{k_z}, & x^{i_y}y^{j}, & x^{i_z}z^{k-1} & \text{if} & T=\TUR,\\
   x^{i} y^{j_y}, & x^{i} z^{k_z+1}, & x^{i_y}y^{j-1}, & x^{i_z}z^{k} & \quad & T=\TUL.
 \end{array}\]

\end{lemma}

\begin{proof}
The proof is a matter of straightforward computation. It follows directly by considering each case separately
% (cf. Figure~\ref{figure:2 valley jig} and
cf.~Lemma\eqref{directions of G-isgaw transformations}).
\end{proof}
Note that the $G-$igsaw transformation of  a $\GG$ with two valleys
may have only one valley.

%\begin{example}\label{example:2 valleys for 1/28(1,5,23)}
%Consider the case of $\frac{1}{28}(1,5,23)$ singularity.
%The $\GG$ $\Gamma=\spann (x^4y,x^4z,x^2y^4,xz^3)$ has two valleys:
%$y-$valley equal to $x^2y$ and $z-$valley equal to $xz.$ The right
%transformation is obtained by setting $u=y^5,v=x^2z$ in the
%definition of the $G-$igsaw transformation
%(Definition~\eqref{definition of Gigsaw}). This means that the
%monomials $x^2z,x^3z,x^4z$ in $\Gamma$ are replaced by
%$y^5,xy^5,x^2y^5,$ respectively. The upper transformation of $\GG$
%is spanned by $x^4y,x^2y^5,xz^3$ and it has $y-$valley equal to
%$x^2y$ and $z-$valley equal to $x.$
%\end{example}

%\begin{figure}[h]\label{figure:2 valley jig} \centering \includegraphics[width=0.8\textwidth]{jigsaw4.eps}
%\caption{Transformations of the $2-$valley $\GG$ for
%$\frac{1}{28}(1,5,23).$} \label{figure:jigsaw}
%\end{figure}

\begin{corollary}\label{transformations of 2 valley into 2 valley remarks}
Let $\Gamma$ be a $\GG$ spanned by $x^i y^{j_y},x^i z^{k_z},x^{i_y}y^j,x^{i_z}z^k$ with 2 valleys: $y-$valley equal to
$v=x^{i_y}y^{j_y}$ and $z-$valley equal to $w=x^{i_z}z^{k_z}.$ If $j_y,k_z\ge 1$ then
\[\begin{array}{ll}
\TR(\TULG)=\Gamma,&\TUL(\TRG)=\Gamma,\\
\TL(\TURG)=\Gamma,&\TUR(\TLG)=\Gamma,
\end{array}\]
that is right and upper left (resp.~left and upper right)
transformations are inverse operations. Moreover,
if $j,k,j-j_y,k-k_z\ge 2$ then
\[\TUL(\TURG)=\TUR(\TULG),\]
that is upper left and upper right transformations commute.
\end{corollary}

\begin{corollary}\label{C:G-igsaws are spanned by}
Let $\Gamma$ be a $\GG$ spanned by $x^i y^{j_y},$ $x^i z^{k_z},$
$x^{i_y}y^j,$ $x^{i_z}z^k,$ with 2 valleys: $y-$valley equal to
$v=x^{i_y}y^{j_y}$ and $z-$valley equal to $w=x^{i_z}z^{k_z}.$ Let
$\Gamma'=\TUR^m(\TUL^n(\Gamma)),$ where $m+n\le
\min\{j,k,j-j_y,k-k_z,\}.$ Then $\Gamma'$ is spanned by $x^{i}
y^{j_y+m},  x^{i} z^{k_z+n},  x^{i_y}y^{j-n}, x^{i_z}z^{k-m}.$ If
$m+n< \min\{j,k,j-j_y,k-k_z\}$ then $\Gamma'$ has two valleys. If
$m+n=\min\{j,k,j-j_y,k-k_z\}$ then $\Gamma'$ has one valley (one of
the monomials $x^{i} y^{j_y+m},  x^{i} z^{k_z+n}, x^{i_y}y^{j-n},
x^{i_z}z^{k-m}$ spanning $\Gamma'$ is redundant).
\end{corollary}

\section{Triangles of transformations and primitive $\GGs$}\label{triangles of transformation}

In this section we introduce primitive $\GGs,$ which have a
particular shape. Every primitive $\GG$ such gives rise to a family
of $\GGs,$ called here a triangle of transformations. It will turn
out that most $\GGs$ belong to some triangle of
transformations. We define a sequence of primitive $\GGs$ containing
every primitive $\GG$ for fixed integers $r$ and $a.$

\begin{definition}
Let $\Gamma$ be a $\GG$ with two valleys, spanned by $x^i y^{j_y},$ $x^i z^{k_z},$
$x^{i_y}y^j,$ $x^{i_z}z^k.$ The set
\[\Theta(\Gamma)=\{\TUR^m(\TUL^n(\Gamma))\;\vert\; m+n\le\min\{j,k,j-j_y,k-k_z\}\}\]
will be called \emph{triangle of transformations} of $\Gamma.$
\end{definition}

The union of the supports of $\GGs$ belonging to the set
$\Theta(\Gamma)$ is a simplicial cone (see Corollary~\eqref{C:same
signs}), hence we call $\Theta(\Gamma)$ a triangle of
transformations.

\begin{definition}\label{definition of primitive G-set}
A $\GG$ $\Gamma$ is called \emph{primitive} if it has a $y-$valley
equal to $x^{i_y}$ and a $z-$valley equal to $x^{i_z}$ for some
nonnegative $i_y,i_z.$
\end{definition}

%\begin{figure}[h] \centering \includegraphics[width=0.2\textwidth]{primitex4.eps}
%\caption{A primitive $\GG$ for $r=14,a=5$ with valleys $1$ and
%$x^3.$} \label{figure:sample primitive}
%\end{figure}

The name primitive is justified by the fact that every $\GG$ with
two valleys belong to a triangle of transformations of some
primitive $\GG.$ This fact will follow from the Main Theorem.

\begin{definition}\label{the primitive G-set}
For fixed coprime integers $r,a$ define let $\Gamma_1$
be a $\GG$ spanned by $x,y^{b-1},z^{r-b-1},$
where $b\in\{1,\ldots,r-1\}$ is as an inverse of $a$ modulo $r.$
\end{definition}

The $\GG$ $\Gamma_1$ is primitive and the monomial $x$ is simultaneously its $y-$valley and $z-$valley.

%\begin{figure}[h]\label{figure:primitive} \centering
%\includegraphics[width=0.6\textwidth]{G1primit4.eps}
%\caption{The $\GG$ $\Gamma_1$ for $r=14,a=5.$}
%\end{figure}

\begin{lemma}\label{L:numbers of GGs in a triangle of transf}
Let $\Gamma$ be a primitive $\GG$ spanned by
$x^i,x^{i_y}y^j,x^{i_z}z^k.$ Then $\Theta(\Gamma)$ consists of
$\binom{\min\{j,k\}+2}{2}$ $\GGs.$
\end{lemma}

\begin{proof}
It is clear from definition of $\Theta(\Gamma).$
\end{proof}

\begin{lemma}\label{L:primitive from primitive}
Let $\Gamma$ be a primitive $\GG$ spanned by
$x^i,x^{i_y}y^j,x^{i_z}z^k.$ Suppose that $j<k$ (resp. $k<j$). The
$\GG$  $\TU(\TUR^j(\Gamma))$ (resp.  $\TU(\TUL^k(\Gamma))$) is
spanned by $x^{i+i_z+1},$ $x^i y^j,$ $x^{i_z}z^{k-(j+1)}$ (resp.
$x^{i+i_y+1},$ $x^i y^{j-(k+1)},$ $x^i z^k$). Moreover, if $j<k-1$
(resp. $k<j-1$) it is primitive.
\end{lemma}

\begin{proof}
Assume that $j<k.$ The $\GG$ $\TUR^j(\Gamma)$ is spanned by $x^i
y^j, x^{i_z}z^{k-j}$ and it has one $z-$valley equal to $x^{i_z}$ by
Lemma~\eqref{transformations of 2 valley into 2 valley}. To finish
the proof apply Lemma~\eqref{classification of upper
transformations} to the $\GG$ $\TUR^j(\Gamma).$
\end{proof}

The preceding lemma allows us to define a sequence of primitive
$\GGs.$

\begin{definition}\label{D:sequence of primitive GG}
If $\Gamma_n$ is a primitive $\GG$ we set:
\[ \Gamma_{n+1}= \left\{
\begin{array}{ccc}
  \TU(\TUR^{j_n}(\Gamma_n)) & \text{if} & j_n<k_n, \\
  \TU(\TUL^{k_n}(\Gamma_n)) & \text{if} & j_n>k_n,
\end{array}
\right.\] where $j_n,k_n$ denote the nonnegative numbers such that
$\Gamma_n$ is spanned by the monomials $x^{i_n},$
$x^{i_{y,n}}y^{j_n},$ $x^{i_{z,n}}z^{k_n}$ for some
$i_n,i_{y,n},i_{z,n}\ge 0.$
\end{definition}

Observe, that if $j_n-k_n=\pm 1$ for some $n$ then $\Gamma_{n+1}$ is
not primitive and the recursion stops.

\begin{corollary}\label{C:formulas}
The numbers $j_n,k_n$ satisfy the following formulas:
\begin{align*}
j_1+1&=b,\\
k_1+1&=r-b,\\
j_{n+1}+1&= \left\{
\begin{array}{lcc}
  j_n+1 & \text{if} & j_n<k_n, \\
  j_n+1-(k_n+1) & \text{if} & j_n>k_n,
\end{array}
\right. \\
k_{n+1}+1&= \left\{
\begin{array}{lcc}
  k_n+1-(j_n+1) & \text{if} & j_n<k_n, \\
  k_n+1 & \text{if} & j_n>k_n.
\end{array}
\right.
\end{align*}
\end{corollary}

Clearly, there is a direct link between the numbers $j_n+1,k_n+1$
and the numbers appearing in the Euclidean algorithm for $b$ and
$r-b.$ This relationship will be exploited later.

\begin{definition}
Let $\Theta(\Gamma)$ be a triangle of transformations of a $\GG$
$\Gamma.$ We define
\[\tTG=\bigcup_{\Gamma'\in\Theta(\Gamma)} \sigma(\Gamma')\]
to be the union of supports of the cones $\sigma(\Gamma'),$ where
$\Gamma'$ runs through the $\GGs$ in $\Theta(\Gamma).$
\end{definition}

To study the location of various cones in the fan $\GHnic\CC^3$ it
is convenient to give names to their rays.

\begin{definition}
Let $\Gamma_n$ be the primitive $\GG$ as defined
in~\eqref{D:sequence of primitive GG}. Denote by $\rho_n$ the common
ray of the cones $\wTheta(\Gamma_n)$ and $\sigma(\Gamma_n).$

Let $\Gamma$ be any $\GG.$ A ray of $\sigma^{\vee}(\Gamma)$ will be
called \emph{upper,(upper) left or right ray} if it dual to the wall
of $\sigma(\Gamma)$ corresponding to the upper,(upper) left or right
transformation, respectively.
\end{definition}

\begin{remark}\label{R:parallel}
Let $\Gamma,\Gamma'$ be any two $\GGs. $ Suppose that the cones
$\sigma(\Gamma)$ and $\sigma(\Gamma')$ intersect either in a
$2-$dimensional face or in a ray. If the cones
$\sigma^{\vee}(\Gamma),\sigma^{\vee}(\Gamma')$ have a common ray
$\rho$ then there exists a $2-$dimensional linear subspace of
$N\otimes \RR$ containing a $2-$dimensional face of $\sigma(\Gamma)$
and of $\sigma(\Gamma')$, both of these dual to the ray $\rho.$
\end{remark}

\begin{lemma}
For any $\GG$ $\Gamma$ with two valleys the set $\tTG$ is a rational
simplicial cone.
\end{lemma}

\begin{proof}
Assume that $\GG$ is spanned by the monomials $x^i y^{j_y},$ $x^i
z^{k_z},$ $x^{i_y}y^j,$ $x^{i_z}z^k$ and let
$l=\min{j,k,j-j_y,k-k_z}.$ Because the upper right and upper left
transformation commute (see Corollary~\eqref{transformations of 2
valley into 2 valley remarks}), by Remark~\eqref{R:parallel} it is
enough to establish
 the three following facts:
\begin{itemize}
\item the right rays of the cones $\sigma^{\vee}(\TUR^n(\Gamma))$ for
$n=0,\ldots,l$ are the same,
\item the left rays of the cones $\sigma^{\vee}(\TUL^n(\Gamma))$ for
$n=0,\ldots,l$ are the same,
\item the upper rays of the cones $\sigma^{\vee}(\TUR^m(\TUL^n(\Gamma)))$ for
$m+n=l$ are the same.
\end{itemize}
These follow from Corollary~\eqref{C:G-igsaws are spanned by}.
\end{proof}

\begin{lemma}\label{L:e_2 or e_3 is a ray}
Let $\Gamma$ be a primitive $\GG$ spanned by
$x^i,x^{i_y}y^j,x^{i_z}z^k.$ If $j<k$ (resp. $k<j$) then $\RR_+ e_2$
(resp. $\RR_+ e_3$) is a ray of $\tTG.$
\end{lemma}

\begin{proof}
Suppose that $j<k.$ The $\GG$ $\Gamma'=\TUL^j(\Gamma)$ is spanned by
the monomials $x^{i} z^{k_z+j},  x^{i_z}z^{k}$ and it has one valley
(see Corollary~\eqref{C:G-igsaws are spanned by}). The upper and
left ray of $\sigma(\Gamma')$ are equal to $x^{i+1}{z^{-k+k_z}}$ and
 ${x^{-i+i_z}}{z^{k+1}},$ respectively. Evidently, the ray of
$\sigma{\vee}(\Gamma'),$ dual to the $2-$dimensional face of
$\sigma^{\vee}(\Gamma)$ spanned by the upper and left ray, is equal
to $\RR_+ e_2.$
\end{proof}

Note, that the cone $\wTheta(\Gamma_i)$ has, besides the ray common
with $\sigma(\Gamma_i)$, two other rays: one equal to either $e_2$
or $e_3$ and the second which belongs to $\sigma(\Gamma_{i+1}).$ We
will investigate how the cones
$\wTheta(\Gamma_i),\wTheta(\Gamma_{i+1})$ fit together depending on
the sign of $(j_i-k_i)(j_{i+1}-k_{i+1}).$

%\begin{figure}[h] \centering \includegraphics[width=0.7\textwidth]{new_triangle.eps}
%\caption{Cone $\tTG$ cut with the hyperplane $e_2^*+e_3^*=28,$ where
%$\Gamma=\spann (x^4,x^2y^5,xz^4)$ for $\frac{1}{28}(1,5,23).$}
%\label{figure:triangle}
%\end{figure}

%\begin{example}
%We consider the case of $\frac{1}{28}(1,5,23)$ singularity. Let
%$\Gamma$ be a primitive $\GG$ spanned by the monomials
%$x^4,x^2y^5,xz^4$
% (cf. Example~(\ref{example:2
%valleys for 1/28(1,5,23)})). By Corollary~(\ref{C:G-igsaws are
%spanned by}), the $\GG$ $\TUR^{4}(\Gamma)$ is spanned by
%$x^4y^4,x^2y^5$ and $\RR_+ e_3$ is a ray of the cone $\tTG.$
%\end{example}

\begin{corollary}\label{C:same signs}
Let $\Gamma_n$ and $\Gamma_{n+1}$ be two primitive $\GGs.$ If
$(j_n-k_n)(j_{n+1}-k_{n+1})>0$ then the union of the supports of the
cones $\widetilde\Theta(\Gamma_n),\widetilde\Theta(\Gamma_{n+1})$ is
a rational simplicial cone.
\end{corollary}

\begin{lemma}\label{L:fitting together}
Let $\Gamma_n$ and $\Gamma_{n+1}$ be two primitive $\GGs.$  Then
$\wTheta(\Gamma_n)\cup\wTheta(\Gamma_{n+1})$ is equal to the cone
spanned by $\rho_n,e_2,e_3$ minus (set-theoretical) the cone spanned
by $\rho_{n+1},e_2,e_3.$
\end{lemma}

\begin{proof}
If $(j_n-k_n)(j_{n+1}-k_{n+1})>0$ this follows from
Corollary~\eqref{C:same signs}. Otherwise, the cones
$\widetilde\Theta(\Gamma_n),\widetilde\Theta(\Gamma_{n+1})$ have a
common ray and a $2-$dimensional face of
$\widetilde\Theta(\Gamma_{n+1})$ is contained in a $2-$dimensional
face of $\widetilde\Theta(\Gamma_n).$ To finish, note that $e_2$ and
$e_3$ generate rays of $\wTheta(\Gamma_n)$ and
$\wTheta(\Gamma_{n+1})$ (up to the order).
\end{proof}

Recall that $\Gyz{l}=\spann (y^{l},z^{r-l-1}).$ We will prove that
the cones $\sigma(\Gyz{l}$ fit nicely together with the cones
$\wTheta(\Gamma_j)$ into the fan of $\GHnic\CC^3.$

\begin{lemma}
The upper transformations of $\Gyz{b-1}$ and $\Gyz{b}$ coincide,
where $b\in\{1,\ldots,r-1\}$ is an inverse of $a$ modulo $r.$ In
fact, they are equal to $\Gamma_1.$
\end{lemma}
\begin{proof}
By definition, the upper transformation of $\Gyz{b-1}$ and $\Gyz{b}$
replaces the monomial $z^{r-b}$ and $y^b$ with the monomial $x,$
respectively.
\end{proof}

\begin{lemma}\label{L:rays of yz G-sets}
The upper rays of the  cones $\sigma^{\vee}(\Gyz{0}),\ldots
\sigma^{\vee}(\Gyz{b-1})$ (resp. $\sigma^{\vee}(\Gyz{b}),\ldots
\sigma^{\vee}(\Gyz{r-1})$) are equal. The $1-$dimensional cone $\RRp
e_1$ is a ray of each the cones $\sigma(\Gyz{i}),$ for
$i=0,\ldots,r-1.$
\end{lemma}

\begin{proof}
The upper ray of the cones $\sigma^{\vee}(\Gyz{0}),\ldots
\sigma^{\vee}(\Gyz{b-1})$ is spanned by $xz^{-r+b}$ and the upper
ray of $\sigma^{\vee}(\Gyz{b}),\ldots \sigma^{\vee}(\Gyz{r-1})$ is
spanned by $xy^b.$ The right and left rays of
$\sigma^{\vee}(\Gyz{l}$ are equal to
$y^{-l}z^{r-l},y^{l+1}z^{r-l-1},$ therefore $\RR_+ e_1$ is a ray of
$\sigma(\Gyz{l}).$
\end{proof}

\begin{corollary}\label{C:sum of yz-cones}
The sets \[\bigcup_{l=0}^{b-1}  \sigma(\Gyz{l}),\quad
\bigcup_{l=b}^{r-1}  \sigma(\Gyz{l})\] are rational cones in
$N\otimes\RR$ spanned by $e_1,e_2,\rho_1$ and $e_1,e_3,\rho_1,$
respectively.
\end{corollary}

\begin{proof}
This follows from Remark~\eqref{R:parallel} and Lemma~\eqref{L:rays
of yz G-sets}.
\end{proof}

\section{Main Theorem and the Euclidean
algorithm}\label{main algorithm}

By Theorem~(\ref{ZNakamury}), when $\Gamma$ varies through all
$\GGs,$ the cones $\sigma(\Gamma)$ form a fan supported on the cone
spanned by $e_1,e_2,e_3.$ Therefore, it is enough to find $\GGs$
different from the $\GG$ $\Gyz{l}$ which does not belong to any
triangle of transformation. By looking at the supports of triangle
transformations, it will turn out that those missing $\GGs$ are
exactly the upper transformations of the last $\GG$ $\Gamma_n$
defined in~\eqref{D:sequence of primitive GG}). With the the help of
the Euclidean algorithm we will be able to give a formula for a
total number of $\GG$ for fixed $r$ and $a.$

\begin{definition}\label{D:last m+1 not primitive}
Let $m$ be an integer such that $\Gamma_{m+1}$ is not primitive
(i.e. $\Gamma_m$ is the last primitive $\GG$ in the sequence defined
in~\eqref{D:sequence of primitive GG}).
\end{definition}

\begin{theorem}[Main Theorem]\label{algorithm}
Let $r,a$ be coprime natural numbers and let $b$ be an inverse of
$a$ modulo $r.$ Let $G$ be a cyclic group of order $r,$ acting on
$\CC^3$ with weights $1,a,r-a.$

If $\Gamma_1,\ldots,\Gamma_{m+1}$ is the sequence from
Definition~\eqref{D:sequence of primitive GG} (that is $\Gamma_n$ is
a primitive $\GG$ unless $n=m+1$) and if $\Gyz{l}=\spann
(y^{r-l-1},z^{l})$ then every $\GG$ either
\begin{itemize}
\item belongs to a
triangle of transformation of some $\Gamma_n$ for $n\le m,$ or
\item is equal to a $\GG$
$\Gyz{l}$ for some $l=1,\ldots,n,$ or
\item is equal to an iterated upper transformation of the $\GG$ $\Gamma_{m+1}.$
\end{itemize}

\end{theorem}

\begin{proof}
The proof uses Nakamura's Theorem~\eqref{ZNakamury}, which asserts
that the union of the supports of the cones $\sigma(\Gamma)$ is equal
to the positive octant in $N\otimes \RR.$ Lemma~\eqref{L:fitting
together} and Corollary~\eqref{C:sum of yz-cones} combined imply
that if a $\GG$ $\Gamma$ neither belongs to some triangle of
transformation nor is equal to $\Gyz{l}$ for some $l$ then the cone
$\sigma(\Gamma)$ is supported in the cone spanned by
$e_2,e_3,\rho_{m+1}.$ On the other hand, the $\GG$ $\Gamma_{m+1}$ is
equal either to $\spann (x^{i_{m+1}},x^{i_{y_{m+1}}}y^{j_{m+1}})$ or
to $\spann (x^{i_{m+1}},x^{i_{z,m+1}}z^{k_{m+1}}),$ cf.
Lemma~\eqref{L:primitive from primitive}. Therefore the $j_{m+1}-$th
or $k_{m+1}-$th iterated upper transformation of $\Gamma_{m+1}$ is
equal to $\Gx=\spann (x^{r-1}).$ Moreover, the $\GGs$ $\TU^l
(\Gamma_{m+1})$ and $\TU^{l+1}(\Gamma_{m+1})$ satisfy assumptions of
the Remark~\eqref{R:parallel}. This shows that the set
\[\bigcup_{l=0}^{\max\{j_{m+1},k_{m+1}\}} \sigma(\TU^l(\Gamma_{m+1}))\]
is a cone generated by $e_2,e_3,\rho_{m+1}$ which concludes the proof.
\end{proof}

\begin{xremark}
The above theorem can be restated in a form of an algorithm
computing the fan of the $\GHnic \CC^3$ for fixed $a$ and $r$
(recall that the $\GHnic \CC^3$ is normal, cf.
Corollary~\eqref{G-Hilb is normal 2}).
\end{xremark}

\begin{xremark}
    The two-stage construction of the $\GHnic\CC^3$ for abelian subgroups in $\SL(3,\CC)$ by Craw and Reid in~\cite{CrawReid:GHSL3} appears to provide a coarse subdivision of the fan of the $\GHnic\CC^3$
    for the subgroup $G$ in $\GL(3,\CC)$ of type $\frac{1}{r}(1,a,r-a).$ The coarse subdivision (i.e. with all interior lines of all triangles of
    transformations removed) is provided by the continued fraction expansions.
\end{xremark}

\begin{lemma}\label{lemma:equalities in eucl algoithm}
Let $p_l,q_l$ be the data of the Euclidean algorithm for the
nonnegative integer numbers $p_1,p_2$ with $\GCD(p_1,p_2)=p_{n+1},$
that is
\[p_i  =  q_i p_{i+1} + p_{i+2}, \quad 0<p_{i+2}<p_{i+1},\]
where $p_{n+1} \neq 0$ and $p_{n+2}  = 0.$

Then \begin{align*}
\sum_{l=1}^{n} q_l p_{l+1}&= p_1 + p_2 - p_{n+1},\\
 \sum_{l=1}^{n} q_l p_{l+1}^2&= p_1 p_2 .
\end{align*}
\end{lemma}

\begin{theorem}\label{the number of G-sets}
Fix some coprime numbers $r$ and $a.$ Let $N$ denote the number of different $\GGs$ for the action of type $\frac{1}{r}(1,a,r-a).$ Then
\[N=\frac{1}{2}(3r+b(r-b)-1).\]
\end{theorem}

\begin{proof}
Denote $\Gamma_l=\spann (x^{i_l}, x^{i_{y,l}}y^{j_l},
x^{i_{z,l}}z^{k_l}).$ The triangle of transformations of $\Gamma_l$
consist of ${\min\{j_l+1,k_l+1\}+1}\choose{2}$ cones (see
Lemma~\eqref{L:numbers of GGs in a triangle of transf}). Therefore
\[ N=r+\max\{j_{m+1}+1,k_{m+1}+1\}+\sum_{l=1}^m {{\min\{j_l+1,k_l+1\}+1}\choose{2}},\]
where the first two terms come from the $\GGs$ $\Gyz{l}$ and the
consecutive upper transformations of $\Gamma_{m+1}.$

Suppose that $b<r-b.$ Let the $p_l$ and $q_l$ be the data of the
Euclidean algorithm for the coprime numbers $p_1=k_1+1=r-b,
p_2=j_1+1=b$ as in Lemma~(\ref{lemma:equalities in eucl algoithm}).
Set $q_0=1.$ In this notation, by the formulas from
Corollary~\eqref{C:formulas},
\begin{gather*}
\min\{j_C+1,k_C+1\}=p_D \\
 \text{for}\quad q_0+\ldots
q_D\le C <q_0+\ldots q_{D+1}.
\end{gather*}
Note that $p_{n+1}=1$ and $q_n=\max\{j_{m+1}+1,k_{m+1}+1\},$ thus
$N=r+q_n p_{n+1}+\frac{1}{2}\sum_{l=1}^{n-1} (q_l p_{l+1}^2+q_l
p_{l+1}).$ This, by simple computation, implies the assertion.

\end{proof}

%\begin{xremark}\label{number of smooth cones}
%By similar computations the number of cones $\sigma(\Gamma)$
%generated by three independent vectors in the fan of $\GHnic$ scheme
%is equal to $2r-1.$ This suggests a link between the Danilov
%resolution and the $\GH$ scheme.
%\end{xremark}

\section{Example}\label{section:example GHilbert}
\begin{figure}[ht!] \centering \includegraphics[width=0.9\textwidth]{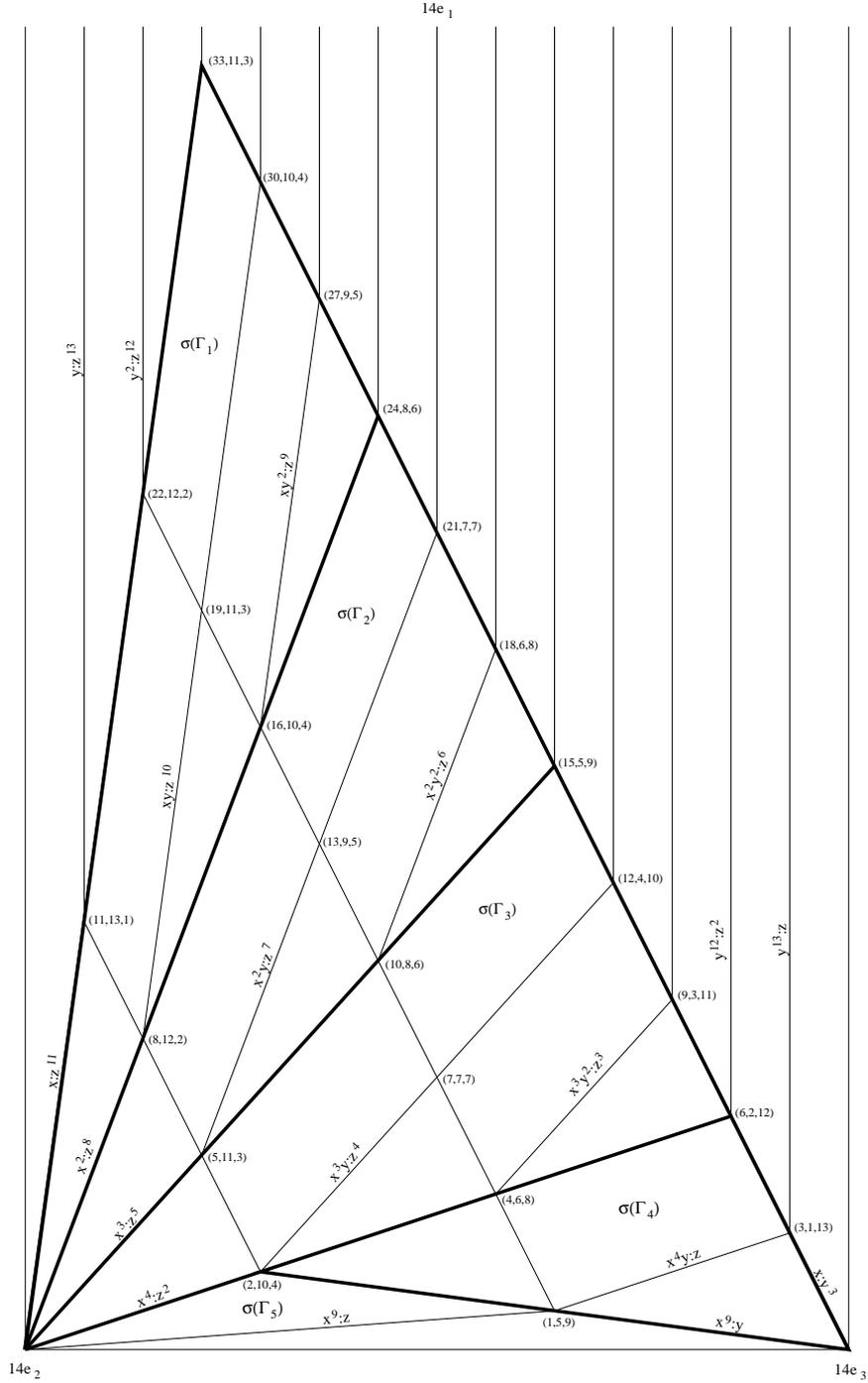}
\caption{The fan of $\GH$ scheme for $r=14,a=5$ intersected with
hyperplane $e^{*}_2+e^{*}_3=14.$} \label{figure:a=5}
\end{figure}
By Theorem~(\ref{algorithm}), for $a=5,r=14$ every $G-$set, different from $\Gyz{i},$ belongs to a triangle
of transformation of the primitive $G-$sets
\begin{align*}
\Gamma_1&=\spann (x,y^{2},z^{10}), \\
\Gamma_2&=\spann (x^2,xy^{2},z^{7}), \\
\Gamma_3&=\spann (x^3,x^2y^{2},z^{4}), \\
\Gamma_4&=\spann (x^4,x^3y^2,z), \\
%\Gamma_5&=\spann (x^8,x^4z), \\
\end{align*}
or is an upper transformation of
\[\TU(\Gamma_5)=\Gx. \]
There are $37$  different $\GGs.$
Figure~\ref{figure:a=5} shows the fan of $\GH,$ where $e_1$ the
ray generated by $e_1$ is drawn at "infinity". The ratios along
lines denote rays of the corresponding cones $\sigma^{\vee}(\Gamma)$
(up to an inverse in the multiplicative notation). Triangles of
transformations are marked with thick line.

\bibliography{literatura}
\bibliographystyle{alpha}
\nocite{Hartshorne:geometry,Fulton:toric,GSVerdier,Reid:Lacoresspondance,Kedzierski:CohGH}
\end{document}